\documentclass[12pt]{report}
\usepackage{amssymb,amsmath,makeidx,verbatim}
\newcommand{\R}{\mathbb{R}}

\newcommand{\C}{\mathbb{C}}
\newcommand{\Z}{\mathbb{Z}}

\newcommand{\be}{\begin{enumerate}}
\newcommand{\ee}{\end{enumerate}}
\newcommand{\bq}{\begin{eqnarray*}}
\newcommand{\eq}{\end{eqnarray*}}

\begin{document}
\newcommand{\disp}{\displaystyle}
\thispagestyle{empty}
\begin{center}
\textsc{Fourier Transform of Schwartz Algebras on Groups in the Harish-Chandra class\\}
\ \\
\textsc{Olufemi O. Oyadare}\\
\ \\
Department of Mathematics,\\
Obafemi Awolowo University,\\
Ile-Ife, $220005,$ NIGERIA.\\
\text{E-mail: \textit{femi\_oya@yahoo.com}}\\
\end{center}
\begin{quote}
{\bf Abstract.} {\it It is well-known that the Harish-Chandra transform, $f\mapsto\mathcal{H}f,$ is a topological isomorphism
of the spherical (Schwartz) convolution algebra $\mathcal{C}^{p}(G//K)$ (where $K$ is a maximal compact subgroup of any arbitrarily
chosen group $G$ in the Harish-Chandra class and $0<p\leq2$) onto the (Schwartz) multiplication algebra $\bar{\mathcal{Z}}({\mathfrak{F}}^{\epsilon})$ (of $\mathfrak{w}-$invariant members of $\mathcal{Z}({\mathfrak{F}}^{\epsilon}),$ with $\epsilon=(2/p)-1$). The same cannot however be said of the full Schwartz convolution algebra $\mathcal{C}^{p}(G),$ except for few specific examples of groups (notably $G=SL(2,\R)$) and for some notable values of $p$ (with restrictions on $G$ and/or on $\mathcal{C}^{p}(G)$). Nevertheless the full Harish-Chandra Plancherel inversion formula on $G$ is known for all of $\mathcal{C}^{2}(G)=:\mathcal{C}(G).$ In order to then understand the structure of Harish-Chandra transform more clearly and to compute the image of $\mathcal{C}^{p}(G)$ under it (without any restriction) we derive an absolutely convergent series expansion (in terms of known functions) for the Harish-Chandra transform by an application of the full Plancherel inversion formula on $G.$ This leads to a computation of the image of $\mathcal{C}(G)$ under the Harish-Chandra transform which may be seen as a concrete realization of Arthur's result and be easily extended to all of $\mathcal{C}^{p}(G)$ in much the same way as it is known in the work of Trombi and Varadarajan.}
\end{quote}
\ \\
\ \\
\ \\
$\overline{2010\; \textmd{Mathematics}}$ Subject Classification: $43A85, \;\; 22E30, \;\; 22E46$
\ \\
\ \\
{\bf \S1. Introduction.} Let $G$ be a reductive group in the \textit{Harish-Chandra class} where $\mathcal{C}^{p}(G)$ is the \textit{Harish-Chandra-type} Schwartz algebra on $G,$ $0<p\leq2,$ with $\mathcal{C}^{2}(G)=:\mathcal{C}(G).$ The image of $\mathcal{C}^{p}(G)$ under the (Harish-Chandra) \textit{Fourier transform} on $G$ has been a pre-occupation of harmonic analysts since the emergence of Arthur's thesis $[1a],$ where the Fourier image of $\mathcal{C}(G)$ was characterized for connected non-compact semisimple Lie groups of real rank one. Thereafter Eguchi $[3a.]$ removed the restriction of the real rank and considered non-compact real semisimple $G$ with only one conjugacy class of \textit{Cartan subgroups} as well as the Fourier image of $\mathcal{C}^{p}(G/K),\;[3b.],$ while Barker $[2.]$ considered $\mathcal{C}^{p}(SL(2,\R))$ as well as $\mathcal{C}^{0}(SL(2,\R)).$ The complete $p=2$ story for any real reductive $G$ is contained in Arthur $[1b,c].$ The most successful general result along the general case of $p$ is the well-known \textit{Trombi-Varadarajan Theorem} $[9.]$ which characterized the image of $\mathcal{C}^{p}(G//K),$ $0<p\leq2,$ for a maximal compact subgroup $K$ of a connected semisimple Lie group $G$ as a (Schwartz) multiplication algebra $\bar{\mathcal{Z}}({\mathfrak{F}}^{\epsilon})$ (of $\mathfrak{w}-$invariant members of $\mathcal{Z}({\mathfrak{F}}^{\epsilon}),$ with $\epsilon=(2/p)-1$); thus subsuming the works of Ehrepreis and Mautner $[4.]$ and Helgason $[6.].$ However the characterization of the image of $\mathcal{C}^{p}(G)$ for reductive groups $G$ in the Harish-Chandra class has not yet been achieved due to failure of the method (successfully employed in $[1b,c]$ and $[10.]$) of generalizing from the real rank one case. Nevertheless the \textit{Plancherel inversion formula} is already known for $\mathcal{C}(G),$ where $G$ is a reductive group in the Harish-Chandra class, $[10.].$

This paper contains a derivation of an absolutely convergent series expansion (in terms of known functions) for the \textit{Harish-Chandra transform} and its application in constructing an explicit image of $\mathcal{C}(G)$ under it, for reductive groups, $G,$ in the Harish-Chandra class. These results give explicit realization of the \textit{abstract} computations of Arthur $[1.]$ and show the direct contributions of the Plancherel inversion formula to the understanding of Harish-Chandra transform. As a Corollary we deduce a more explicit form of the $p=2$ Trombi-Varadarajan Theorem.
\ \\
\ \\
{\bf \S2. Preliminaries.} Let $G$ be a group in the Harish-Chandra class. That is $G$ is a locally compact group with the properties that $G$ is reductive, with Lie algebra $\mathfrak{g},$ $[G:G^{0}]<\infty,$ where $G^{0}$ is the connected component of $G$ containing the identity, in which the analytic subgroup, $G_{1},$ of $G$ defined by $\mathfrak{g}_{1}=[\mathfrak{g},\mathfrak{g}]$ is closed in $G$ and of finite center and in which, if $G_{\C}$ is the adjoint group of $\mathfrak{g}_{\C},$ then $Ad(G)\subset G_{\C}.$ Such a group $G$ is endowed with a \textit{Cartan involution,} $\theta,$ whose fixed points form a \textit{maximal} compact subgroup, $K,$ of $G$ $[5.].$ $K$ meets all connected components of $G,$ in particular $G^{0}.$

We denote the \textit{universal enveloping algebra} of $\mathfrak{g}_{\C}$ by $\mathcal{U}(\mathfrak{g}_{\C}),$ whose members may be viewed either as left or right invariant differential operators on $G.$ We shall write $f(x;a)$ for the left action $(af)(x)$ and $f(a;x)$ for the right action $(fa)(x)$ of $\mathcal{U}(\mathfrak{g}_{\C})$ on functions $f$ on $G.$ Let $\mathcal{C}(G)$ represents the space of $C^{\infty}-$functions on $G$ for which $$\sup_{x \in G}\mid f(b;x;a) \mid \Xi^{-1}(x)(1+\sigma(x))^{r}<\infty,$$ for $a,b \in \mathcal{U}(\mathfrak{g}_{\C})$ and $r>0.$ Here $\Xi$ and $\sigma$ are well-known \textit{elementary} spherical functions on $G.$ $\mathcal{C}(G)$ is known to be a Schwartz algebra under convolution while $\mathcal{C}(G//K),$ consisting of the spherical members of $\mathcal{C}(G),$ is a closed commutative subalgebra.

Let $\hat{G}$ represent the set of equivalence classes of \textit{irreducible unitary representations} of $G.$ If $G_{1}$ is non-compact then the support of the Plancherel measure does not exhaust $\hat{G}.$ We write $\hat{G_{t}}$ for this support, which generally contains a discrete part, $\hat{G_{d}}$ ($\neq\emptyset,$ if $rank(G) = rank(K)$), and a continuous part, $\hat{G_{t}}\setminus \hat{G_{d}}$ ($\neq\emptyset,$ always). $G$ has finitely many conjugacy classes of \textit{Cartan subgroups,} which may be represented by the set $\{A_{1},\cdots,A_{r}\}.$

We can write $$A_{i}=A_{i,I}\cdot A_{i,\R},$$ where $A_{i,I}=A\cap K$ is a maximal compact subgroup of $A_{i}$ and $A_{i,\R}$ (whose Lie algebra shall be denoted as $\mathfrak{a}_{i,\R}$) is a vector subgroup with $\theta(a)=a^{-1},$ $a \in A_{i,\R}.$ There are parabolic subgroups $P_{i}$ of $G$ whose \textit{Langlands decompositions} are of the form $P_{i}=M_{i}A_{i,\R}N_{i}.$ Each of the subgroups $M_{i}$ satisfies all the requirements on $G.$

The \textit{full} Harish-Chandra Plancherel inversion formula on $G$ states that there are uniquely defined non-negative continuous functions, $C_{i},$ on $(\hat{M_{i}})_{d}\times \mathfrak{a}_{i,\R}$ such that each $C_{i}$ is of at most polynomial growth in which $C_{i}(s\omega:s\nu)=C_{i}(\omega:\nu),$ $s \in W(G,A_{i})$ (the Weyl group of the pair $(G,A_{i})$) and, for $f \in \mathcal{C}(G),$ we have $$f(1)=\sum^{r}_{i=1}(\sum_{\omega \in (\hat{M_{i}})_{d}}\int_{\mathfrak{a}^{*}_{i,\R}}C_{i}(\omega:\nu)\Theta_{\omega,\nu}(f)d\nu),$$ with absolute convergence, $[10.]$. Here $\Theta_{\omega,\nu}$ represents the \textit{distributional characters} of irreducible unitary representations, $\pi_{\omega,\nu},$ on $G.$

More explicitly, the support of the \textit{reduced dual} $\widehat{G}_{r}$ of $G$ split as $\widehat{G}^{d}_{r}\bigsqcup\widehat{G}^{M}_{r},$ it is populated only by the \textit{discrete,} $\pi_{\omega},\;\omega\in\widehat{G}^{d}_{r}$ and \textit{principal series,} $\pi_{(\sigma,\lambda)},\;(\sigma,\lambda)\in\widehat{M}\times
\R=\widehat{G}^{M}_{r}$ of representations of $G$ and there exist a number $\beta(\omega),\;\omega\in\widehat{G}^{d}_{r}$ (called the \textit{formal degree} of $\omega$) and a non-negative function $\beta:\widehat{M}\times\R\rightarrow[0,\infty):(\sigma,\lambda)\mapsto\beta(\sigma,\lambda)$ such that, for any $f\in\mathcal{C}(G),$ the split version of the full Harish-Chandra Plancherel inversion formula on $G$ is given as $$f(1)=\sum_{\omega\in\widehat{G}^{d}_{r}}\beta(\omega)\Theta_{\omega}(f)+
\sum_{\sigma\in\widehat{M}}\int^{\infty}_{0}\beta(\sigma,\lambda)\Theta_{\sigma,\lambda}(f)d\lambda.$$

The first form of the Plancherel inversion formula given above (which bears a striking resemblance with members of the space $L^{2}_{\wp}(\hat{G})$ given in $[1b,c]$) may be re-written as $\;\;\;\;f(1)=\sum^{r}_{i=1}(\sum_{\omega \in (\hat{M_{i}})_{d}}\int_{\mathfrak{a}^{*}_{i,\R}}C_{i}(\omega:\nu)\Theta_{\omega,\nu}(f)d\nu)$ $$=\sum^{r}_{i=1}(\sum_{\omega \in (\hat{M_{i}})_{d}}\int_{\mathfrak{a}^{*}_{i,\R}}\| a_{\wp_{i}}(\omega,\nu) \|^{2} \mu_{\omega}(\nu)d\nu),$$ where the measurable functions $$(\omega,\nu)\mapsto a_{\wp_{i}}(\omega,\nu)$$  (of Arthur $[1b,c]$) on $(\hat{M_{i}})_{d}\times \mathfrak{a}_{i,\R}$ are realizable as those for which $$\| a_{\wp_{i}}(\omega,\nu) \|^{2} \mu_{\omega}(\nu)=C_{i}(\omega:\nu)\Theta_{\omega,\nu}(f)\;\mbox{and}\;\|a_{\wp}\|^{2}=f(1),$$ for any $f \in \mathcal{C}(G)$ and $\mu_{\omega}-$almost everywhere. This gives a first suggestion that the Plancherel inversion formula may contribute to the realization of the image of $\mathcal{C}(G).$ It then means that Arthur's map $(\omega,\nu)\mapsto a_{\wp_{i}}(\omega,\nu)$ are functions whose $L^{2}-$norm satisfies the above prescriptions of being realizable in terms of members of $\mathcal{C}(G).$ This remark may prove useful in the eventual realization of these functions. It will however be clear from our main Theorem that an explicit expression for $a_{\wp_{i}}(\omega,\nu)$ (though of interest in its own right) would not be necessary in the computation of $\mathcal{H}(\mathcal{C}(G)).$

We shall denote the (elementary) spherical functions on $G$ by $\varphi_{\lambda},$ $\lambda \in \mathfrak{a}^{*}_{\C}.$ We know that the  \textit{functional equation} $$\varphi_{s\lambda}(x)=\varphi_{\lambda}(x)$$ holds ($[10.],\;p.\;365$) for all $s \in \mathfrak{w}=\mathfrak{w}(\mathfrak{g}_{\C},\mathfrak{a}_{\C}),\;x \in G.$ Since $$\mathcal{C}(G)*\mathcal{C}(G)\rightarrow \mathcal{C}(G)$$ is a continuous map, then the function $x\mapsto s_{\lambda,f}(x),$ $x \in G,$ given as $$s_{\lambda,f}(x):=(f*\varphi_{\lambda})(x)$$ is a Schwartz function on $G.$ It is also clear that $s_{\lambda,f}(1)=(\mathcal{H}f)(\lambda),$ $[8.]$ and $[9.].$
\ \\
\ \\
{\bf \S3. Harish-Chandra Fourier Transform on $\mathcal{C}(G)$.} We shall start by introducing the following space of functions on $\mathfrak{a}^{*}.$

{\bf Definition 3.1} \textit{Let $ \{A_{1},\cdots,A_{r}\}$ be a complete set of representatives in the conjugacy classes of Cartan subgroups of $G$ and let $\mathcal{C}_{G}(\mathfrak{a}^{*})$ represents the space of functions on $\mathfrak{a}^{*}$ that are topologically spanned by maps of the form $$\lambda\mapsto g_{i,f}(\lambda):=\int_{\mathfrak{a}^{*}_{i,\R}}C_{i}(\omega:\nu)\Theta_{\omega,\nu}(f*\varphi_{\lambda})d\nu,$$ where $i=1,\cdots,r,$ and $f \in \mathcal{C}(G).$}

It is clear (from the fact that $\hat{G_{t}}\setminus \hat{G_{d}}\neq\emptyset$) that $\mathcal{C}_{G}(\mathfrak{a}^{*})\neq\emptyset$ for any reductive group, $G,$ in the Harish-Chandra class. We endow $\mathcal{C}_{G}(\mathfrak{a}^{*})$ with the basic operations of a function space which convert it into a multiplication algebra. Since $f \in \mathcal{C}(G)$ and $\varphi_{\lambda} \in \mathcal{C}(G//K),$ the algebra $\mathcal{C}_{G}(\mathfrak{a}^{*}),$ consisting of absolutely convergent integrals whose decay estimates follow from the above Plancherel inversion formula, is therefore a Schwartz algebra on $\mathfrak{a}^{*}.$

{\bf Lemma 3.2} \textit{$\mathcal{C}_{G}(\mathfrak{a}^{*})$ is invariant under the action of the Weyl group $\mathfrak{w}=\mathfrak{w}(\mathfrak{g}_{\C},\mathfrak{a}_{\C})$ defined as $(s\cdot g)(\lambda)=g(s^{-1}\lambda),$ $s \in \mathfrak{w},$ $g \in \mathcal{C}_{G}(\mathfrak{a}^{*}).$}\\
{\bf Proof.} Note that $(s\cdot g_{i,f})(\lambda)=g_{i,f}(s^{-1}\lambda)=\int_{\mathfrak{a}^{*}_{i,\R}}C_{i}(\omega:\nu)\Theta_{\omega,\nu}(f*\varphi_{s^{-1}\lambda})d\nu
$ $$=\int_{\mathfrak{a}^{*}_{i,\R}}C_{i}(\omega:\nu)\Theta_{\omega,\nu}(f*\varphi_{s_{1}\lambda})d\nu$$ $$\mbox{(since there exists $s_{1} \in \mathfrak{w}$ such that $s_{1}\cdot s=1$ and so that $s^{-1}=s_{1}$)}$$ $$=\int_{\mathfrak{a}^{*}_{i,\R}}C_{i}
(\omega:\nu)\Theta_{\omega,\nu}(f*\varphi_{\lambda})d\nu
=g_{i,f}(\lambda).\;\Box$$

It may also be shown that, if $g \in \mathcal{C}_{G}(\mathfrak{a}^{*})$ and $\lambda_{1},\lambda_{2} \in \mathfrak{a}^{*}$ are such that $s\lambda_{1}=\lambda_{2},$ for some $s \in \mathfrak{w},$ then $g(\lambda_{1})=g(\lambda_{2}).$ The above Lemma implies that $\mathcal{C}_{G}(\mathfrak{a}^{*})^{\mathfrak{w}}=\mathcal{C}_{G}(\mathfrak{a}^{*})$ and that we do not have to consider the $\mathfrak{w}-$invariant subspace of $\mathcal{C}_{G}(\mathfrak{a}^{*}).$ Indeed $\mathcal{C}_{G}(\mathfrak{a}^{*})$ is the $\mathfrak{w}-$invariant subspace of a larger space of functions on $\mathfrak{a}^{*}$ spanned by maps of the form $$\lambda\mapsto g(\lambda):=\int_{\mathfrak{a}^{*}_{i,\R}}C_{i}(\omega:\nu)\Theta_{\omega,\nu}(h_{\lambda})d\nu,$$ for \textit{arbitrary chosen} $h_{\lambda} \in \mathcal{C}(G),$ $i=1,\cdots,r,$ indexed by $\lambda \in \mathfrak{a}^{*}.$

\textbf{Lemma 3.3} \textit{The Harish-Chandra Fourier transform $f\mapsto \mathcal{H}f$ is a map of $\mathcal{C}(G)$ into $\mathcal{C}_{G}(\mathfrak{a}^{*})$}\\
{\bf Proof.} As $f \in \mathcal{C}(G)$ we have that $f*\varphi_{\lambda} \in \mathcal{C}(G)$ for $\lambda \in \mathfrak{a}^{*}.$ Hence by the above Plancherel inversion formula we conclude that $$(f*\varphi_{\lambda})(1)=\sum^{r}_{i=1}(\sum_{\omega \in (\hat{M_{i}})_{d}}\int_{\mathfrak{a}^{*}_{i,\R}}C_{i}(\omega:\nu)\Theta_{\omega,\nu}(f*\varphi_{\lambda})d\nu).$$ That is, $$(\mathcal{H}f)(\lambda)=\sum^{r}_{i=1}(\sum_{\omega \in (\hat{M_{i}})_{d}}\int_{\mathfrak{a}^{*}_{i,\R}}C_{i}(\omega:\nu)\Theta_{\omega,\nu}(f*\varphi_{\lambda})d\nu),$$ which show that $\mathcal{H}:\mathcal{C}(G)\rightarrow \mathcal{C}_{G}(\mathfrak{a}^{*}),$ with all its known properties still intact$.\Box$

A more importatnt realization of the algebra $\mathcal{C}_{G}(\mathfrak{a}^{*})$ is derived as follows.

\textbf{Lemma 3.4} \textit{The algebra $\mathcal{C}_{G}(\mathfrak{a}^{*})$ is the topological linear span of maps of the form $$\lambda\mapsto\int_{\mathfrak{a}^{*}_{i,\R}}C_{i}(\omega:\nu)\pi_{\omega,\nu,\mathfrak{d}_{i}}(f)d\nu$$ where $f\in\mathcal{C}(G)$ and $\mathfrak{d}_{i}$ is the singleton set consisting of the trivial representation of each $M_{i}$ corresponding to the elementary spherical function $\varphi_{\lambda}.$}\\
{\bf Proof.} We only need to show that $$\Theta_{\omega,\nu}(f*\varphi_{\lambda})=\pi_{\omega,\nu,F_{i}}(f).$$ Indeed, we have that $$\Theta_{\omega,\nu}(f*\varphi_{\lambda})=\int_{G}\Theta_{\omega,\nu}(x)(f\ast\varphi_{\lambda})(x)dx$$ $$=\int_{G}\Theta_{\omega,\nu}(x)(\int_{G}f(y)\varphi_{\lambda}(y^{-1}x)dy)dx$$ $$=\int_{G}f(y)(\int_{G}\Theta_{\omega,\nu}(x)\varphi_{\lambda}(y^{-1}x)dx)dy$$ $$
\int_{G}f(y)(\Theta_{\omega,\nu}\ast\varphi_{\lambda})(y)dy=
\int_{G}f(y)\Phi_{\pi_{\omega,\nu},\mathfrak{d}}(y)dy\\
=:\pi_{\omega,\nu,\mathfrak{d}}(f)$$ (in the notation of $[5.],\;p.\;22),$ since the function $\Phi_{\pi_{\omega,\nu},\mathfrak{d}}$ given above as $\Phi_{\pi_{\omega,\nu},\mathfrak{d}}=\Theta_{\omega,\nu}\ast\varphi_{\lambda}$ is a spherical function of type $\mathfrak{d}=\{1\}$ which is associated with $\pi_{\omega,\nu}.\Box$

\textbf{Corollary 3.5} \textit{The Harish-Chandra Fourier transform $f\mapsto \mathcal{H}f$ of $\mathcal{C}(G)$ into $\mathcal{C}_{G}(\mathfrak{a}^{*})$ is given by an absolutely convergent series as $$(\mathcal{H}f)(\lambda)=\sum^{r}_{i=1}(\sum_{\omega \in (\hat{M_{i}})_{d}}\int_{\mathfrak{a}^{*}_{i,\R}}C_{i}(\omega:\nu)\pi_{\omega,\nu,\mathfrak{d}_{i}}(f)d\nu).\Box$$}

The above function $\Phi_{\pi_{\omega,\nu},\mathfrak{d}}$ is worthy of an independent study. It is a spherical function (being a convolution of two spherical functions) and is continuous (due to the continuity of convolution in $\mathcal{C}(G)$). We may normalize it in order to have $\Phi_{\pi_{\omega,\nu},\mathfrak{d}}(1)=1.$ We shall denote the \textit{isotypical subspace} of $\mathcal{C}(G)$ with respect to $\mathfrak{d}=\{1\}$ as $\mathcal{C}_{1}(G).$

\textbf{Lemma 3.6} \textit{There exists $f\in\mathcal{C}_{1}(G)$ such that $\pi_{\omega,\nu,\mathfrak{d}}(f)=1.$}\\
\textbf{Proof.} Since $\pi_{\omega,\nu,\mathfrak{d}}$ is an irreducible representation on $G$ in which the isotypical subspace $\mathcal{C}_{1}(G)$ satisfies $dim(\mathcal{C}_{1}(G))<\infty\;([11],\;p.\;114),$ our result follows from using Lemma $1.3.2$ of $[5].\;p.\;22.\Box$

\textbf{Theorem 3.7} \textit{Arthur's Fourier image $\mathcal{C}(\widehat{G})$ of $\mathcal{C}(G)$ coincides with the algebra $\mathcal{C}_{G}(\mathfrak{a}^{*}).$}\\
\textbf{Proof.} The representation $\pi_{\omega,\nu,\mathfrak{d}}$ satisfies all the requirements in Lemma $1.3.3$ of $[5],\;p.\;22.$ Hence, the topological linear span of the functions $$\lambda\mapsto\int_{\mathfrak{a}^{*}_{i,\R}}C_{i}(\omega:\nu)\pi_{\omega,\nu,\mathfrak{d}_{i}}(f)d\nu$$ (in Lemma $3.4$ above), which is $\mathcal{C}_{G}(\mathfrak{a}^{*}),$ gives the same algebra as that of $\pi_{\omega,\nu}(f),\;f\in\mathcal{C}(G),$ which we already know to be $\mathcal{C}(\widehat{G}).\Box$

We have therefore proved the main Theorem of this paper.

{\bf Theorem 3.8} \textit{The Harish-Chandra Fourier transform, $f\mapsto \mathcal{H}f,$ is a topological algebra isomorphism of $\mathcal{C}(G)$ onto the multiplication algebra $\mathcal{C}_{G}(\mathfrak{a}^{*}).$}\\
\textbf{Proof.} Combine Corollary $3.5$ with Theorem $3.7.\Box$

\textbf{Corollary 3.9} \textit{The map $I_{\mathfrak{d}}:\mathcal{C}(\widehat{G})\rightarrow\mathcal{C}_{G}(\mathfrak{a}^{*})$ defined as $$I_{\mathfrak{d}}(\pi_{\omega,\nu}(f))=\sum^{r}_{i=1}(\sum_{\omega \in (\hat{M_{i}})_{d}}\int_{\mathfrak{a}^{*}_{i,\R}}C_{i}(\omega:\nu)\pi_{\omega,\nu,\mathfrak{d}_{i}}(f)d\nu)$$ is a topological isomorphism$.\Box$}

An immediate consequence of our results is the following which gives an explicit realization of the celebrated Trombi-Varadarajan Theorem when $p=2.$

{\bf Corollary 3.10 (\textmd{p=2} Trombi-Varadarajan Theorem)} \textit{Let ${\mathfrak{F}}^{\epsilon}$ be the Trombi-Varadarajan tubular region in $\mathfrak{a}^{*}$ containing $\mathfrak{a}^{*}_{i,\R},$ where $\epsilon=(2/p)-1$ with $0<p\leq2.$ Then the Trombi-Varadarajan image, $\bar{\mathcal{Z}}({\mathfrak{F}}^{0}),$ is a closed subalgebra of $\mathcal{C}_{G}(\mathfrak{a}^{*})$ and $$\bar{\mathcal{Z}}({\mathfrak{F}}^{0})=\mathcal{C}_{G//K}(\mathfrak{a}^{*}).$$}
{\bf Proof.} Since $\mathcal{C}(G//K)$ is a closed subalgebra of $\mathcal{C}(G)$ and the map $\mathcal{C}(G//K)*\mathcal{C}(G//K)\rightarrow \mathcal{C}(G//K)$ is continuous, Theorem $3.8$ implies that we only need to show that $\mathcal{H}f \in \mathcal{C}_{G}(\mathfrak{a}^{*}),$ for $f \in \mathcal{C}(G//K).$ Now let $f \in \mathcal{C}(G//K)$ and $\lambda \in \mathfrak{a}^{*},$ then $f*\varphi_{\lambda} \in \mathcal{C}(G//K),$ so that the map $$\lambda\mapsto g(\lambda):=\int_{\mathfrak{a}^{*}_{i,\R}}C_{i}(\omega:\nu)\Theta_{\omega,\nu}(f*\varphi_{\lambda})d\nu,$$ is well-defined on $\mathfrak{a}^{*}.$ Hence $g \in \mathcal{C}_{G}(\mathfrak{a}^{*})$ and $\mathcal{H}f \in \mathcal{C}_{G}(\mathfrak{a}^{*}).\;\Box$

The proof of the above Theorem $3.8$ contains a new \textit{formula-representation} of the Harish-Chandra transform as an absolutely convergent series expansion in terms of the Harish-Chandra $c-$functions and global distributions. \textit{This formula shows (for the first time) the direct dependence of the transform on the uniquely defined non-negative continuous functions, $C_{i},$ and characters, $\Theta_{\omega,\nu},$ of the irreducible unitary representations on $G.$} This expansion may further expose other properties of the transform on a closer look. The importance of our approach to harmonic analysis of Schwartz algebra on $G$ is that the (\textit{\textit{character} form} of) Plancherel inversion formula of any reductive (Lie) group leads directly to the computation of the Fourier image of this algebra, an image that would be as explicit as the Plancherel inversion formula itself.

An analogous realization of $\mathcal{C}_{G}(\mathfrak{a}^{*})$ may still be sought through an explicit derivation of the combined decay properties of $$C_{i}(\omega:\nu)\Theta_{\omega,\nu}(f*\varphi_{\lambda})=C_{i}(\omega:\nu)\pi_{\omega,\nu,\mathfrak{d}_{i}}(f),$$ for $f \in \mathcal{C}(G),$ especially for specific $G.$ The continuous inclusion $\mathcal{C}^{p}(G)\subset L^{2}(G)$ makes the above ($L^{2}-$) Plancherel inversion formula and absolutely convergent series expansion of the Harish-Chandra transform available for use in extending Theorem $3.8$ from $\mathcal{C}(G)$ to all of $\mathcal{C}^{p}(G)$ (and hence Corollary $3.10$ from $\bar{\mathcal{Z}}({\mathfrak{F}}^{0})$ to all of $\bar{\mathcal{Z}}({\mathfrak{F}}^{\epsilon})$), in exactly the same way that the results on $\mathcal{C}(SL(2,\R)//SO(2))$ was extended to all of $\mathcal{C}^{p}(SL(2,\R)//SO(2))$ in $[9.].$

A close comparison with formula $(15.1)$ of $[2]$ leads to the following.

\textbf{Conjecture 3.11} \textit{The map $I_{\mathfrak{d}}$ in Corollary $3.9$ is an explicit realization of the inverse Fourier transform.}

\ \\
{\bf \S4. The Example of $SL(2,\R).$} The elementary spherical functions on $G=SL(2,\R)$ are given as the matrix coefficients of class $1$ members of the principal series $\pi_{\omega\nu}$ $([5],\;p.\;103)$ and are found to be the \textit{Legendre functions} $$\varphi_{\nu}(a_{t})=P_{i\nu+\frac{1}{2}}(\cosh t)=\frac{1}{2\pi}\int^{2\pi}_{0}(\cosh t+\sinh t\cos u)^{i\nu+\frac{1}{2}}du,$$ $t\in\R\;([6b],\;p.\;406).$ Distributional characters of $G$ are $G-$invariants defined as $$\Theta_{\pi}(f)=tr\pi(f),$$ for all $f$ for which $\pi(f):=\int_{G}f(x)\pi(x)dx$ is of \textit{trace class.} These characters are completely defined on the only two non-conjugate Cartan subgroups of $G=SL(2,\R)$ namely the non-compact type $T=MA$ (with $M=\{\pm I\}$ and $A=\{a_{t}:=\left(\begin{array}{cc} e^{t} & 0 \\ 0 & e^{-t} \end{array}\right):\;t\in\R\}$) and the compact type given as $B:=K=\{k_{\theta}:=\left(\begin{array}{cc} \cos\theta & \sin\theta \\ -\sin\theta & \cos\theta \end{array}\right):\theta\in[0,2\pi)\},$ with corresponding Haar measures of $d(\pm)dt$ and $\frac{1}{2\pi}d\theta,$ where both of functions $(e^{i\theta}-e^{-i\theta})\Theta_{\pi}(\theta)$ and $\pm\mid e^{t}-e^{-t}\mid\Theta_{\pi}(\pm a_{t})$ are locally integrable functions with respect to $(B,d(\pm)dt)$ and $(K,\frac{1}{2\pi}d\theta,)$ respectively.

These facts suggest exponential expressions for both $(e^{i\theta}-e^{-i\theta})\Theta_{\pi}(\theta)$ and $\pm\mid e^{t}-e^{-t}\mid\Theta_{\pi}(\pm a_{t}),$ where $\pi=\pi_{\sigma,\nu}$ is the principal series representations (in which we now write $\Theta_{\pi_{\sigma,\nu}}$ simply as $\Theta_{\sigma,\nu}$) and $\pi=\pi^{\pm}_{n}$ is the principal series representations (in which we now write $\Theta_{\pi^{\pm}_{n}}$ simply as $\Theta_{n}$) of $G.$ Indeed we have that $\widehat{M}=\{\sigma_{\pm}\}$ (with $\sigma_{+}(\pm I)=1$ and $\sigma_{-}(\pm I)=\pm1$) and $\widehat{K}=\Z^{*}=\Z\setminus\{0\},$ so that $$\Theta_{\sigma,\nu}(x)=\left\{\begin{array}{ll} \;\;\;\;\;\;\;\;\;\;\;\;0, & \mbox{if}\;
x\;\;\mbox{is conjugate to}\;\;k_{\theta},\\
\frac{\sigma(\pm)(e^{\nu(\log a_{t})}+e^{-\nu(\log a_{t})})}{\mid e^{t}-e^{-t} \mid}, & \mbox{if}\;
x\;\;\mbox{is conjugate to}\;\;\pm a_{t},
\end{array}\right.$$ and
$$\Theta_{n}(x)=\left\{\begin{array}{ll} \;\;\;\;\;\;\frac{\alpha e^{-\alpha(n-1)i\theta}}{(e^{i\theta}-e^{-i\theta})}, & \mbox{if}\;
x\;\;\mbox{is conjugate to}\;\;k_{\theta},\\
(\pm)^{n}\frac{e^{(n-1)t}(1-sgn\;t)+e^{-(n-1)t}(1+sgn\;t)}{2\mid e^{t}-e^{-t} \mid}, & \mbox{if}\;
x\;\;\mbox{is conjugate to}\;\;\pm a_{t},
\end{array}\right.$$ (Here $\alpha=1$ when $\pi=\pi^{-}_{n}$ and $\alpha=-1$ when $\pi=\pi^{+}_{n}$) $[7],\;p.\;342,\;343.$

With $f\in\mathcal{C}(SL(2,\R))$ we have the Plancherel inversion formula given as $$f(1)=\frac{1}{2\pi}\sum_{k\in\Z^{*}}\mid k\mid\Theta_{k}(f)+ (\frac{1}{2\pi})^{2}
\sum_{\sigma\in\widehat{M}}\int_{Re\nu=0}\beta(\sigma,\lambda)\Theta_{\sigma,\lambda}(f)d\nu.$$ with $$\beta(\sigma,\nu)=\left\{\begin{array}{ll} (\frac{\lambda\pi i}{2})\tan(\frac{\lambda\pi}{2}), & \mbox{if}\;
\sigma=\sigma_{+},\\
(-\frac{\lambda\pi i}{2})\cot(\frac{\lambda\pi}{2}), & \mbox{if}\;
\sigma=\sigma_{-},
\end{array}\right.$$ Hence, from Corollary $3.5,$ we have $$(\mathcal{H}f)(\lambda)=\frac{1}{2\pi}\sum_{k\in\Z^{*}}\mid k\mid\pi_{k,1}(f)+ (\frac{1}{2\pi})^{2}
\sum_{\sigma\in\widehat{M}}\int_{Re\nu=0}\beta(\sigma,\nu)\pi_{\sigma,\nu,1} d\nu,$$ which by Lemma $3.3,$ becomes $$(\mathcal{H}f)(\lambda)=\frac{1}{2\pi}\sum_{k\in\Z^{*}}\mid k\mid\Theta_{k}(f*\varphi_{\lambda})+ (\frac{1}{2\pi})^{2}
\sum_{\sigma\in\widehat{M}}\int_{Re\nu=0}\beta(\sigma,\nu)\Theta_{\sigma,\nu}(f*\varphi_{\lambda})d\nu.$$ Hence
$$(\mathcal{H}f)(\lambda)\leq\mid\frac{1}{2\pi}\sum_{k\in\Z^{*}}\mid k\mid\Theta_{k}(f*\varphi_{\lambda})+ (\frac{1}{2\pi})^{2}
\sum_{\sigma\in\widehat{M}}\int_{Re\nu=0}\beta(\sigma,\nu)\Theta_{\sigma,\nu}(f*\varphi_{\lambda})d\nu\mid$$
$$\leq\frac{1}{2\pi}\sum_{k\in\Z^{*}}\mid k\mid\cdot\mid\Theta_{k}(f*\varphi_{\lambda})\mid+ (\frac{1}{2\pi})^{2}
\sum_{\sigma\in\widehat{M}}\int_{Re\nu=0}\mid\beta(\sigma,\nu)\mid\cdot\mid\Theta_{\sigma,\nu}(f*\varphi_{\lambda})\mid d\nu.$$

By a calculation similar to that of $[11],\;p.\;271,$ there exist continuous seminorms, $\lambda\mapsto v_{1}(\lambda):=v_{1}(\varphi_{\lambda})$ and $\lambda\mapsto v_{2}(\lambda):=v_{2}(\varphi_{\lambda})$ such that $$\mid\Theta_{\sigma,\nu}(f*\varphi_{\lambda})\mid\leq (1+\mu^{2})^{-2}\cdot v_{1}(\varphi_{\lambda})$$ and $$\mid\Theta_{k}(f*\varphi_{\lambda})\mid\leq m^{-4}\cdot v_{2}(\varphi_{\lambda}).$$ Thus $$(\mathcal{H}f)(\lambda)\leq v_{2}(\lambda)\cdot\sum_{m\neq0}\mid m\mid^{-3}+ cv_{1}(\lambda)\cdot\int_{0}^{\infty}\mu(1+\mu^{2})^{-2}d\mu,$$ proving that $\mathcal{H}f$ is a continuous seminorm on $\R,$ as expected in Theorem $3.8.$ A similar computation suggests Conjecture $3.11.$\\
\ \\
{\bf   References.}
\begin{description}
\item [{[1.]}] Arthur, J. G., $(a.)$ \textit{Harmonic analysis of tempered distributions on semisimple Lie groups of real rank one}, Ph.D. Dissertation, Yale University, $1970;$ $(b.)$ \textit{Harmonic analysis of the Schwartz space of a reductive Lie group I,} mimeographed note, Yale University, Mathematics Department, New Haven, Conn; $(c.)$ \textit{Harmonic analysis of the Schwartz space of a reductive Lie group II,} mimeographed note, Yale University, Mathematics Department, New Haven, Conn.
        \item [{[2.]}] Barker, W. H.,  $L^p$ harmonic analysis on $SL(2, \R),$ \textit{Memoirs of American Mathematical Society,} $vol.$ {\bf 76} , no.: {\bf 393}. $1988$
            \item [{[3.]}] Eguchi, M., $(a.)$ The Fourier Transform of the Schwartz space on a semisimple Lie group, \textit{Hiroshima Math. J.}, $vol.$ \textbf{4}, ($1974$), pp. $133-209.$ $(b.)$ Asymptotic expansions of Eisenstein integrals and Fourier transforms on symmetric spaces, \textit{J. Funct. Anal.} \textbf{34}, pp. $167 - 216.$
                \item [{[4.]}] Ehrenpreis, L. and Mautner, F. I., Some properties of the Fourier transform on semisimple Lie groups, I,
                \textit{Ann. Math.}, $vol.$ $\textbf{61}$ ($1955$), pp. 406-439; II, \textit{Trans. Amer. Math. Soc.}, $vol.$ $\textbf{84}$ ($1957$), pp. $1-55;$ III, \textit{Trans. Amer. Math. Soc.}, $vol.$ $\textbf{90}$ ($1959$), pp. $431-484.$
                    \item [{[5.]}] Gangolli, R. and Varadarajan, V. S., \textit{Harmonic analysis of spherical functions on real reductive groups,} Ergebnisse der Mathematik und iher Genzgebiete, $vol.$ {\bf 101}, Springer-Verlag, Berlin-Heidelberg. $1988.$
                    \item [{[6.]}] Helgason, S., $(a)$  A duality for symmetric spaces with applications to group representations, \textit{Advances in Mathematics,} $vol.$ $\textbf{5}$ ($1970$), pp. $1-154.$ $(b)$ \textit{Differential geometry and symmetric spaces,} Academic Press, New York, $1962.$
                        \item [{[7.]}] Knapp, A.W., \textit{Representation theory of semisimple groups; An overview based on examples,} Princeton University Press, Princeton, New Jersey. $1986.$
                \item [{[8.]}] Oyadare, O. O., On harmonic analysis of spherical convolutions on semisimple Lie groups, \textit{Theoretical Mathematics and Applications}, $vol.$ $\textbf{5},$ no.: {\bf 3}. ($2015$), pp. $19-36.$
                \item [{[9.]}] Trombi, P. C. and Varadarajan, V. S., Spherical transforms on semisimple Lie groups, \textit{Ann. Math.,} $vol.$ {\bf 94}. ($1971$), pp. $246$-$303.$
                        \item [{[10.]}] Varadarajan, V. S., Eigenfunction expansions on semisimple Lie groups, in \textit{Harmonic Analysis and Group Representation}, (A. Fig$\grave{a}$  Talamanca (ed.)) (Lectures given at the $1980$ Summer School of the \textit{Centro Internazionale Matematico Estivo (CIME)} Cortona (Arezzo), Italy, June $24$ - July $9.$ vol. \textbf{82}) Springer-Verlag, Berlin-Heidelberg. $2010,$ pp. $351-422.$
                            \item [{[11.]}] Varadarajan, V. S., \textit{An introduction to harmonic analysis on semisimple Lie groups,} Cambridge Studies in Advanced Mathematics, \textbf{161},  Cambridge University Press, $1989.$
                        \end{description}
\end{document}